\documentclass[12pt,twoside]{amsart}
\usepackage{amscd,amsfonts,amsmath,amssymb,amsthm,latexsym}
\usepackage{enumerate,array,srcltx,amsmath}
\usepackage{leqno}
\usepackage[dvips]{graphicx}
\usepackage{multicol}

\usepackage{caption}
\usepackage{subcaption}

\usepackage{color}

\definecolor{darkblue}{rgb}{0,0,0.545098}
\definecolor{darkgreen}{rgb}{0,0.392157,0}

 \newtheorem{thm}{Theorem}
 
 \newtheorem{lem}[thm]{Lemma}
 \newtheorem{prop}[thm]{Proposition}
 \theoremstyle{definition}
 
 \theoremstyle{remark}
 \newtheorem{rem}{Remark}

 \newenvironment{pf}{\vspace{0.5em}\noindent\textbf{Proof.} }{\quad \hfill $\Box$ \\ \vspace{0.5em}\\}

\def\cH{\mathcal{H}}
\def\ln{\lambda^n}
\def\Re{\textrm{Re}}
\def\Im{\textrm{Im}}
\def\limRe{-\dfrac{1}{\beta}}
\def\limCo{-\dfrac{1}{2}\left(\dfrac{1}{\alpha}-\dfrac{1}{\beta}\right)}

\def\wwd{\widetilde{d}(\mu_n)}

\def\R{\mathbb{R}}
\def\C{\mathbb{C}}

\title[Optimal scalar products in the Standard Linear Viscoelastic Model]{Optimal scalar products in the \\ Standard Linear Viscoelastic Model$^*$}
\thanks{$^*$This work is partially supported by grants MTM2011-27739-C04-01 (Spain) and MTM2011-27739-C04-03 (Spain). The authors are part of the Catalan research group 2014 SGR 1083.}
\author[M. Pellicer] {M. Pellicer$^{\dag}$}
\thanks{$^{\dag}$Dpt. d'Inform\`atica, Matem\`atica Aplicada i Estad\'{\i}stica (EPS), Universitat de Girona, Campus de Montilivi, 17071 Girona, Catalunya (Spain), e-mail: martap@imae.udg.edu}
\author[J. Sol\`{a}-Morales] {J. Sol\`{a}-Morales$^{\ddag}$}
\thanks{$^{\ddag}$Departament de
Matem\`atica Aplicada 1, Universitat Polit\`ecnica de Catalunya,
Av. Diagonal 647, 08028 Barcelona, Catalunya (Spain), e-mail:
jc.sola-morales@upc.edu}

\date{February 23, 2015}

\begin{document}

\maketitle
\begin{abstract}
We study the third order in time linear dissipative wave equation known as the Standard Linear Viscoelastic Model, that appears also as the
linearization of the so-called Moore-Gibson-Thompson equation in Nonlinear Acoustics. We complete the description in \cite{Trig2013} of the spectrum of the generator of
the corresponding group of operators and show that, apart from some exceptional values of the parameters,
this generator can be made to be a normal operator with a new scalar product, with a complete set of orthogonal eigenfunctions. Using this property we also obtain sharper decay estimates for the solutions
as $t\to\infty$, both when the operator is normal or not.

\end{abstract}

\section{Introduction and statement of results}

In this paper we study the third order in time dissipative abstract wave equation
\begin{equation}\label{eq}
(u+\alpha u_{t})_{tt}+L (u+ \beta u_{t}) = 0, \ \ \textrm{with } \alpha,\beta>0
\end{equation}
where $L$ is a self-adjoint, strictly positive operator in a Hilbert space $H$  with compact resolvent. Recall that in this situation the eigenvalues $\mu_n$ of $L$ are positive, increasing, tending to $\infty$, semi-simple and the corresponding eigenfunctions $\phi_n$ are an orthonormal family.  A typical case is when $L=-a^2\Delta$ and $\Delta$ is Laplace's operator with Dirichlet boundary conditions in a bounded domain $\Omega$ with a regular boundary, and in that case $H=L^2(\Omega)$.

The dissipative case corresponds to $\alpha/\beta<1$. It is known that in this case all solutions tend exponentially to zero (see \cite{Lasiecka2011}, \cite{Trig2013}). For $\alpha/\beta\ge 1$ it can be shown, by looking for the appropriate eigenvalues and eigenfunctions or by seeing that the energy is non-decreasing (see formula (\ref{energy}) below, or \cite{Gorain2010}, \cite{MR2013}), that there exist solutions that do not tend to zero when $t\to\infty$.

When $L=-a^2\Delta$, this equation is called the Standard Linear Viscoelastic model, and represents the linear deformations of a viscoelastic solid with an approach that is considered to be more realistic than the usual Kelvin-Voigt model. See \cite{Gorain2010} and \cite{MR2013} and the references therein for a discussion of this model. The same equation appears as a linearization for a classical model in nonlinear acoustics, where it is known under the name of Moore-Gibson-Thompson equation. See \cite{Lasiecka2011} and \cite{Trig2013} and their references. Other recent studies on this equation can be seen in \cite{Andrade} and \cite{Kalantarov}, for example. The paper \cite{Conejero} deals with the case $\alpha/\beta > 1$.

The second-order in time strongly damped wave equation, that would correspond to taking $\alpha=0$ in (\ref{eq}), was studied by the authors in \cite{JDE} some years ago, motivated by our previous work \cite{JMAA} on a viscoelastic model. In that paper we proved that the infinitesimal generator of the semigroup associated to that equation was self-adjoint in a particular new scalar product, that was shown explicitly, provided that the dissipation coefficient was large enough (overdamping regime). Several optimal decay estimates for the solutions were then obtained by using that new scalar product.

The goal of the present paper is also to show the existence of a suitable scalar product associated to the semigroup defined by (\ref{eq}) for some ranges of the parameters. In this scalar product the infinitesimal generator of the semigroup turns out to be a {\em normal} operator, hence admitting an orthonormal basis of eigenfunctions. The necessary and sufficient condition for the existence of this scalar product, and the normality of the operator, is given in Theorem \ref{thm1}, to be stated below. This fact will allow us to obtain the optimal exponential decay rate of the solutions. The main point that will have to be checked in the construction of this scalar product is that it defines a new norm that is equivalent to the natural one. This will be done in several functional spaces where the equation (\ref{eq}) defines a semigroup (in fact, a group). The second main result of the present paper is given in Theorem \ref{thm2}. It states the optimal decay rate for the solutions, and not only the growth bound of the semigroup, both when the scalar product considered in Theorem \ref{thm1} exists or not.

Writing $U=(u,v,w)^T=(u,u_t,u_{tt})^T$ we can write (\ref{eq}) as the following first order evolution equation:

\begin{equation}\label{eq:evolution}
  \dfrac{dU}{dt}= A U\ ,\ U\in D(A), \textrm{   with    }  A U=\left(\begin{array}{c}
    v \\ w \\ -\dfrac{1}{\alpha} L (u+\beta v)-\dfrac{1}{\alpha} w
  \end{array}\right)
\end{equation}

The operator $A$ can be defined in several functional spaces. According to \cite{MR2013}, if $L=-a^2\Delta$ with Dirichlet boundary conditions then the first possibility is that:

\begin{equation}\label{espaisMR}
\begin{array}{lcl}
 \cH &=& H^1_0(\Omega)\times H^1_0(\Omega)\times L^2(\Omega)  \\
 D(A) &=& \{(u,v,w),w\in H^1_0(\Omega),u+\beta v \in H^2(\Omega)\cap H^1_0(\Omega)\}
\end{array}
\end{equation}

It can be seen that, in this case, $(A,D(A))$ defines a $C^0$-semigroup (in fact a group) and that it is dissipative when $\alpha<\beta$. This last part can be seen using the energy functional associated to the following scalar product:
\begin{equation}\label{eq:prodMR}
  \begin{array}{c}
  \left< (u_1,v_1,w_1),(u_2,v_2,w_2)\right> =  \\ \\
  \int_{\Omega} (v_1+\alpha w_1)\overline{(v_2+\alpha w_2)} +
  a^2\int_{\Omega} \nabla(u_1+\alpha v_1)\nabla(\overline{u_2+\alpha v_2}) +
  a^2\alpha(\beta-\alpha)\int_{\Omega} \nabla v_1 \nabla\overline{ v_2 }
  \end{array}
\end{equation}
Observe that if $(u,v,w)\in \cH$ this energy is well defined. From \cite{MR2013} or \cite{Gorain2010} we know that if $(u,v,w)\in D(A)$ then
\begin{equation*}\label{energy}
   \dfrac{d E(t)}{dt}= -a^2(\beta-\alpha)\int_{\Omega} |\nabla v|^2
\end{equation*}

which exhibits the dissipativeness of the operator when $\alpha<\beta$.

But this is not the only possible functional framework. According to \cite{Trig2013}, some possible functional settings are:
\begin{equation}\label{espaisTrig}
\begin{array}{lcl}
  \cH_1 &=& D(L^{1/2})\times D(L^{1/2})\times H  \\
  \cH_2 &=& D(L)\times D(L)\times D(L^{1/2})  \\
  \cH_3 &=& D(L)\times D(L^{1/2})\times H  \\
  \cH_4 &=&D(L^{3/2})\times D(L)\times D(L^{1/2}),
\end{array}
\end{equation}
with the corresponding domains for $A$. Observe that, actually, $\cH_1=\cH$ defined in \eqref{espaisMR} if $L=-a^2\Delta$ with Dirichlet boundary conditions.

The normality of the infinitessimal generator $A$ in a new explicit metric is given in Theorem \ref{thm1}. Associated to \eqref{eq} we define the following numbers $m_1,m_2$, which, as it will be seen in Section \ref{sec:espectre}, are the zeroes of a certain Cardano discriminant, in fact that of the characteristic equation \eqref{eq:char_eq} below:
\begin{equation}\label{m1m2}
m_1 =\alpha\dfrac{-C_1-\sqrt{C_2}}{8\beta^3},\ m_2 =\alpha\dfrac{-C_1+\sqrt{C_2}}{8\beta^3}
\end{equation}
  with
  \begin{equation}\label{eq:AB}
    C_1 = 27-18\left(\dfrac{\beta}{\alpha}\right) -\left(\dfrac{\beta}{\alpha}\right)^2  ,\ \
    C_2 =  C_1^2-64\left(\dfrac{\beta}{\alpha}\right)^3.
  \end{equation}

\begin{thm}\label{thm1}
If $\mu_n\ne m_1,m_2$ for all the eigenvalues $\mu_n$ of $L$, then in each of the spaces $\cH_i$ given in \eqref{espaisTrig} we can define a new equivalent and explicit scalar product $\langle \cdot,\cdot\rangle_{G_i}$ where the operator $A$ becomes a normal operator. Also, there exists a set $\{ \Psi_j^{n,i}, j=1,2,3, n=1,\ldots, \infty\}$ of eigenfunctions of $A$ which is orthonormal in the corresponding new scalar product and complete in $\cH_i$.

Conversely, in the cases where one of the eigenvalues of $L$ coincides with one of the two numbers $m_1$ or $m_2$ (including the case $m_1=m_2$) then the operator $A$ can not be made to be normal in any scalar product.
\end{thm}

\begin{rem}\label{skip}
Associated to each $\mu_n$ eigenvalue of $L$, there exist three corresponding eigenvalues of $A$, named $\lambda^n_1, \lambda^n_2,\lambda^n_3$, the three solutions of the characteristic equation:
\begin{equation}\label{eq:char_eq}
\alpha \lambda^3 + \lambda^2+\beta\mu_n\lambda+\mu_n=0.
\end{equation}
The role of the numbers $m_1,m_2$ will become more clear in Section \ref{sec:espectre} (see Proposition \ref{descriptioneig}) when we show that:
\begin{enumerate}
  \item if $m_1<\mu_n<m_2$, the three of $\lambda^n_1, \lambda^n_2,\lambda^n_3$ are real;
  \item if $\mu_n=m_1$ or $\mu_n=m_2$, they are also real, but two of them are equal and not semi-simple;
  \item if $\mu_n=m_1=m_2$ then $\lambda_1^n=\lambda_2^n=\lambda_3^n\in\mathbb{R}$ with algebraic multiplicity equal to three;
  \item otherwise, $\lambda^n_1$ will be real and $\lambda^n_2=\overline{\lambda^n_3}\in\mathbb{C}\setminus\mathbb{R}$.
\end{enumerate}
 None of the situations 1, 2 or 3 was considered to be possible in \cite{Trig2013}, due to a small error in the analysis of the characteristic equation. This error becomes more important if we consider $\alpha$ variable and near $0$ since in this case the situation $\mu_n=m_1$ or $\mu_n=m_2$ will happen for infinitely many values of $\alpha$.

\end{rem}

\begin{rem}
  In contrast with the strongly damped case $\alpha=0$ studied in \cite{JDE} there is no hope to obtain that the operator is self-adjoint in any new metric, because of the existence of nonreal eigenvalues when $\alpha>0$. So, the property of being a normal operator is the best we can expect.
\end{rem}

\begin{thm}\label{thm2}$\ $
\begin{enumerate}[i)]
\item Suppose that $A$ is a normal operator in the new scalar product $G$ obtained in Theorem \ref{thm1}. Then, any solution $U(t)$ of \eqref{eq:evolution} decays exponentially in the corresponding norm as
\begin{equation*}
  \| U(t)\|_G\leq e^{\sigma_{max} t}\| U(0)\|_G, \textsl{  for }t\geq 0
\end{equation*}
where $\sigma_{max}=\sigma_{max}(A)<0$ is the supremum respect to $n$ (which sometimes is a maximum)  of the real parts of the solutions of the characteristic equation \eqref{eq:char_eq}, that is, the real part of the sometimes called dominant spectrum of $A$ (see Proposition \ref{dominant} below for a description of $\sigma_{max}$).
\item This decay is optimal in the sense that for each $\omega<\sigma_{max}$ there exist solutions $U(t)$ such that
$$\| U(t)\|_G\ e^{-\omega t}\to \infty \textsl{ as } t\to\infty.$$
\item On the other hand, if we are in the situation where, according to Theorem \ref{thm1}, we can not have a new scalar product where the operator $A$ is normal, the previous optimal exponential decay rate result still holds in another suitable norm.
    \end{enumerate}
\end{thm}

\begin{rem}
Theorem \ref{thm2} recalls the decay results of \cite{Trig2013}, but we slightly improve them in the sense that what they prove is that $\sigma_{max}$ is the so-called growth bound of the semigroup, that is
\begin{equation*}
   \inf \{ \omega\in\mathbb{R} ; \| e^{At}\| \leq M_{\omega} e^{\omega t} \ \forall t\geq 0\}
\end{equation*}
(in the usual norm). We show that this infimum is, in fact, a minimum and also that we can take $M_{\sigma_{max}}=1$ in a suitable equivalent norm. This is true even for the case when $\sigma_{max}=-1/\beta$, the essential spectrum of $A$ (see Propositions \ref{essential} and \ref{dominant}).
\end{rem}

The results just stated will be developed and proved in the following sections. More concretely, in Section \ref{sec:normal} we will show the existence and form of the new equivalent scalar product for the given range of the parameters and also prove the normality of the operator in this case and the non-normality otherwise (proof of Theorem \ref{thm1}). Also in Section \ref{sec:normal}, we will prove the optimal exponential decay rate of the solutions in all the cases (that is, Theorem \ref{thm2}). The next Section \ref{sec:espectre} will be devoted to the description of the spectrum of $A$, completing the results of \cite{Trig2013}.

\section{Description of the spectrum of $A$}\label{sec:espectre}

In this section we make an accurate description of the spectrum of $A$, $\sigma(A)$. Most of the results can be found in \cite{Trig2013}, but we include them here for a better global comprehension. Nevertheless, there are some differences, that come mainly from the cases where there are three real eigenvalues associated to the same value of $\mu_n$, that can even be algebraically double or triple, cases that were skipped in \cite{Trig2013}, as it has been said in Remark \ref{skip}. Our description is summarized in the following three propositions.

\begin{prop}[The essential spectrum, see \cite{Trig2013}, Theorem 3.2]\label{essential}
In the four functional settings considered in \eqref{espaisTrig} the essential spectrum of the operator $A$ is
$$\sigma_{ess}(A)=\left\{ -\dfrac{1}{\beta}\right\}.$$
\end{prop}

The definition of essential spectrum can be found in \cite{Henry} or \cite{Gohb-K}.

\begin{prop} [Description of the eigenvalues]\label{descriptioneig}
  The operator $A$ has an infinite number of isolated eigenvalues all of them with finite algebraic multiplicity. More concretely, for each $\mu_n$, $n\in\mathbb{N}$, eigenvalue of $L$, there exist three corresponding eigenvalues of $A$, named $\ln_1$, $\ln_2$ and $\ln_3$, the three solutions of the corresponding characteristic equation \eqref{eq:char_eq}. Moreover, if $\lambda_j^n=a_n+ib_n$ is nonreal, then $a_n$ and $b_n$ satisfy:
        \begin{equation}\label{eq:partreal}
        8\alpha a_n^3+8 a_n^2+2 a_n\left( \dfrac{1}{\alpha}+\beta\mu_n\right)+\mu_n\left(\dfrac{\beta}{\alpha}-1\right)=0,
        \end{equation}
        \begin{equation}\label{eq:partim}
        b_n^2=\dfrac{1}{\alpha} \left( 3\alpha a_n^2 +2 a_n +  \beta\mu_n\right).
        \end{equation}

  Under the dissipativeness condition $0<\alpha<\beta$, these eigenvalues satisfy the following:

  \begin{enumerate}
    \item
    \begin{enumerate}[(a)]
    \item If $\frac{1}{9}<\frac{\alpha}{\beta}<1$, then
    for all $n$ one of the eigenvalues is real and the other two are complex conjugated: $\ln_1\in\mathbb{R}$ and $\ln_2=\overline{\ln_3}\in\mathbb{C}\setminus\mathbb{R}$.

    \item If $0<\frac{\alpha}{\beta}<\frac{1}{9}$, the same happens except, maybe, for a finite number of values of $n$.
    In this case $0<\frac{\alpha}{\beta}<\frac{1}{9}$, the roots of a certain Cardano discriminant $m_1,m_2$ given in \eqref{m1m2} (see the proof for details) are real and satisfy $0<m_1<m_2$. Then, if $m_1<\mu_n<m_2$ (this can happen only for a finite number of values of $n$) then $\ln_1,\ln_2,\ln_3\in\mathbb{R}$. The following exceptional case can also happen (double root case): if  there exist $n_1$ or $n_2$ such that $\mu_{n_1}=m_1$ or $\mu_{n_2}=m_2$, then $\lambda^{n_1}_2=\lambda^{n_1}_3\in\mathbb{R}$ or $\lambda^{n_2}_2=\lambda^{n_2}_3\in\mathbb{R}$ is an eigenvalue with algebraic multiplicity equal to two. Alternatively, if $\mu_n\notin [m_1,m_2]$ then $\ln_1\in\mathbb{R}$ and $\ln_2=\overline{\ln_3}\in\mathbb{C}\setminus\mathbb{R}$, as in (a).

    \item If $\frac{\alpha}{\beta}=\frac{1}{9}$, then $m_1=m_2=\frac{3}{\beta^2}>0$ and the same as in (a) also happens except if there exists $n_1$ such that $\mu_{n_1}=\frac{3}{\beta^2}$. In this exceptional case,  $\lambda^{n_1}_1=\lambda^{n_1}_2=\lambda^{n_1}_3=-\frac{3}{\beta}\in\mathbb{R}$
        is an eigenvalue with algebraic multiplicity equal to three.
    \end{enumerate}

    \item

    \begin{enumerate}[(a)]
        \item If $\lambda$ is a real eigenvalue of $A$, then
        \begin{equation}\label{eq:fita_lambda}
                -\dfrac{1}{\alpha} < \lambda <-\dfrac{1}{\beta}.
              \end{equation}
        If $\lambda$ is nonreal, then
        \begin{equation}\label{eq:fita_part re}
                \limCo<\Re(\lambda)<0.
              \end{equation}

      \item If $\mu_n<\mu_m$ and $\mu_n,\mu_m\notin[m_1,m_2]$, then $\Re(\lambda^n_2)>\Re(\lambda^m_2)$.

      \item The limits of the sequences of eigenvalues are the following:
            \begin{equation}\label{eq:lim_l1}
              \lim_{n\to\infty} \ln_1=\limRe,
            \end{equation}
            which is not an eigenvalue, but the only element in the essential spectrum of $A$ (see Proposition \ref{essential}), 
            \begin{equation}\label{eq:lim_Rel2}
            \lim_{n\to\infty} \Re(\ln_2)=\limCo
            \end{equation}
            and $\lim_{n\to\infty} \Im(\ln_2)=\infty$ with
            \begin{equation}\label{eq:lim_Iml2}
              \Im (\ln_2) = \sqrt{\dfrac{\beta}{\alpha}}\sqrt{\mu_n} + o\left(\sqrt{\mu_n}\right).
            \end{equation}

            \item Also, one has that $\lim_{n\to\infty} \Re(\ln_2)$ is lower than $\lim_{n\to\infty} \ln_1$ (respectively, equal or higher), if $\dfrac{\alpha}{\beta}$ is lower than $\dfrac{1}{3}$ (respectively, equal or higher).

  \end{enumerate}

\end{enumerate}
\end{prop}

The proof of this Proposition \ref{descriptioneig} is presented below. The proof focusses in the cases not considered in \cite{Trig2013}.

In the next proposition we describe the dominant part of the spectrum, that is the part with the highest real part. This real part will be named $\sigma_{max}(A)$ (or simply $\sigma_{max}$ when there is no confusion).

\begin{figure}
        \centering
        \begin{subfigure}[b]{0.4\textwidth}
                \includegraphics[width=\textwidth]{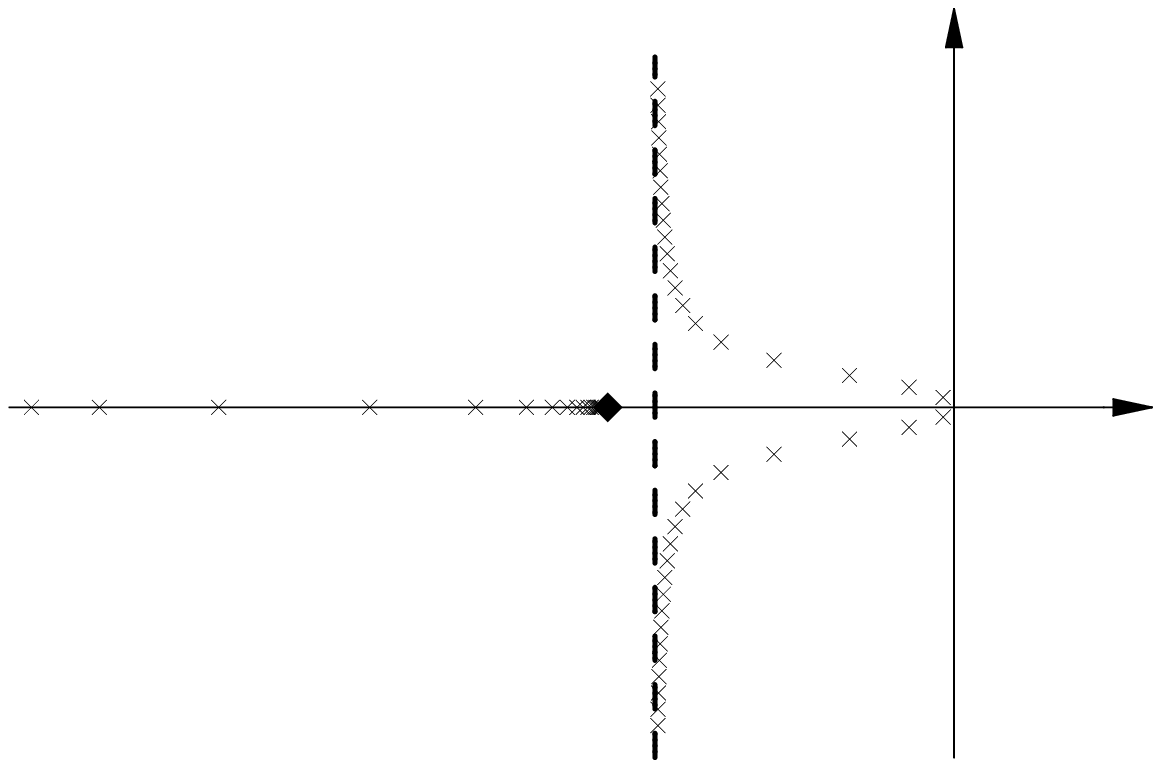}
                \caption{$\sigma_{max}=\Re (\lambda_2^1)$}
                \label{fig:fig1}
        \end{subfigure}\hspace{2cm}
        \begin{subfigure}[b]{0.4\textwidth}
                \includegraphics[width=\textwidth]{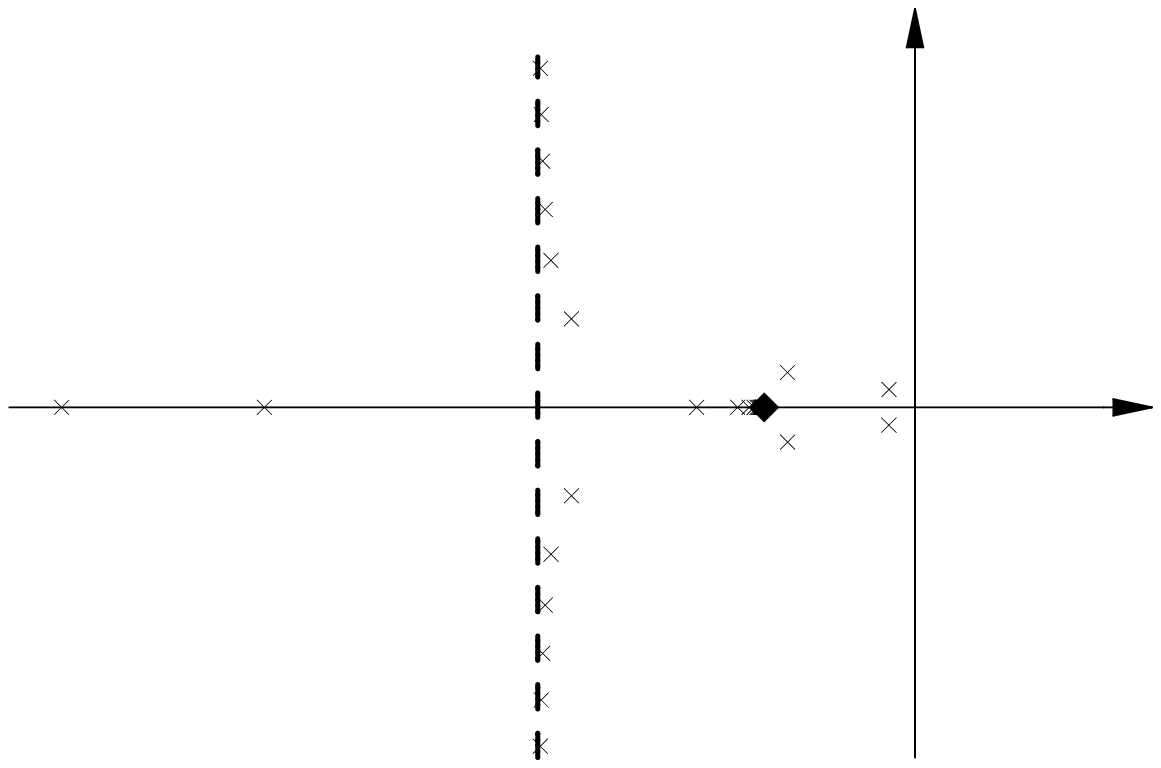}
                \caption{$\sigma_{max}=\Re (\lambda_2^1)$}
                \label{fig:fig2}
        \end{subfigure}\\ \vspace{0.2cm}
        \begin{subfigure}[b]{0.4\textwidth}
                \includegraphics[width=\textwidth]{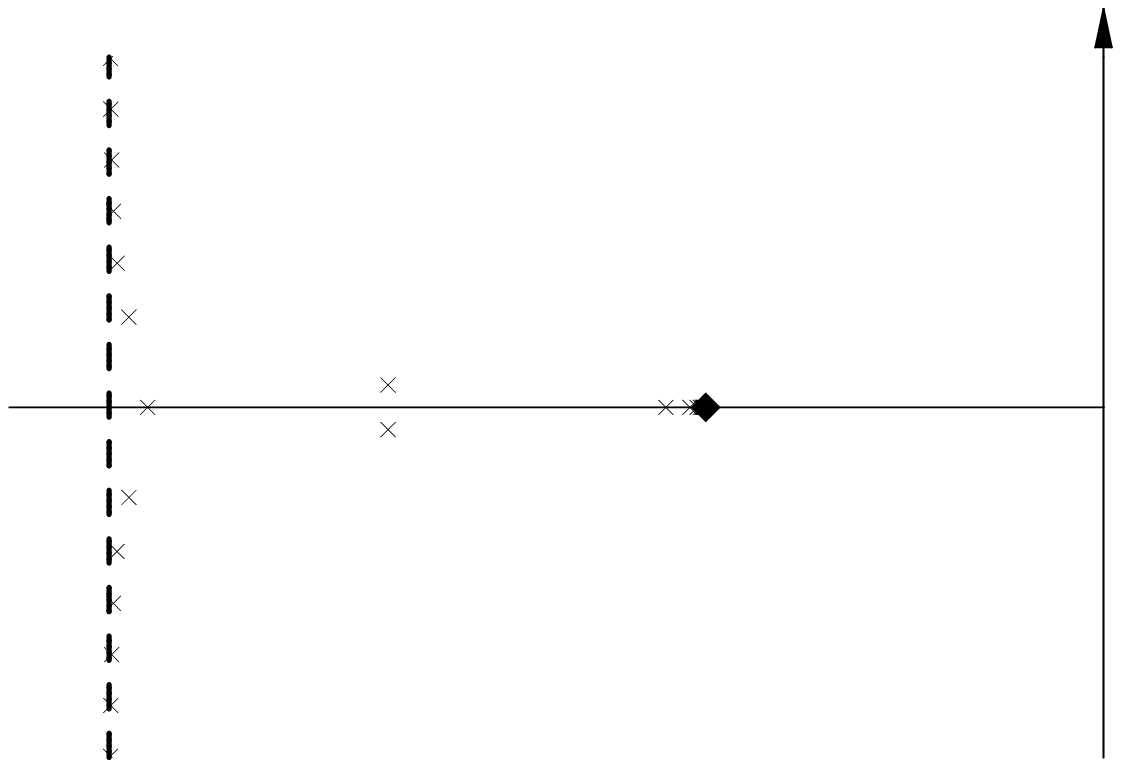}
                \caption{$\sigma_{max}=-\frac{1}{\beta}$}
                \label{fig:fig3}
        \end{subfigure}\hspace{2cm}
        \begin{subfigure}[b]{0.4\textwidth}
                \includegraphics[width=\textwidth]{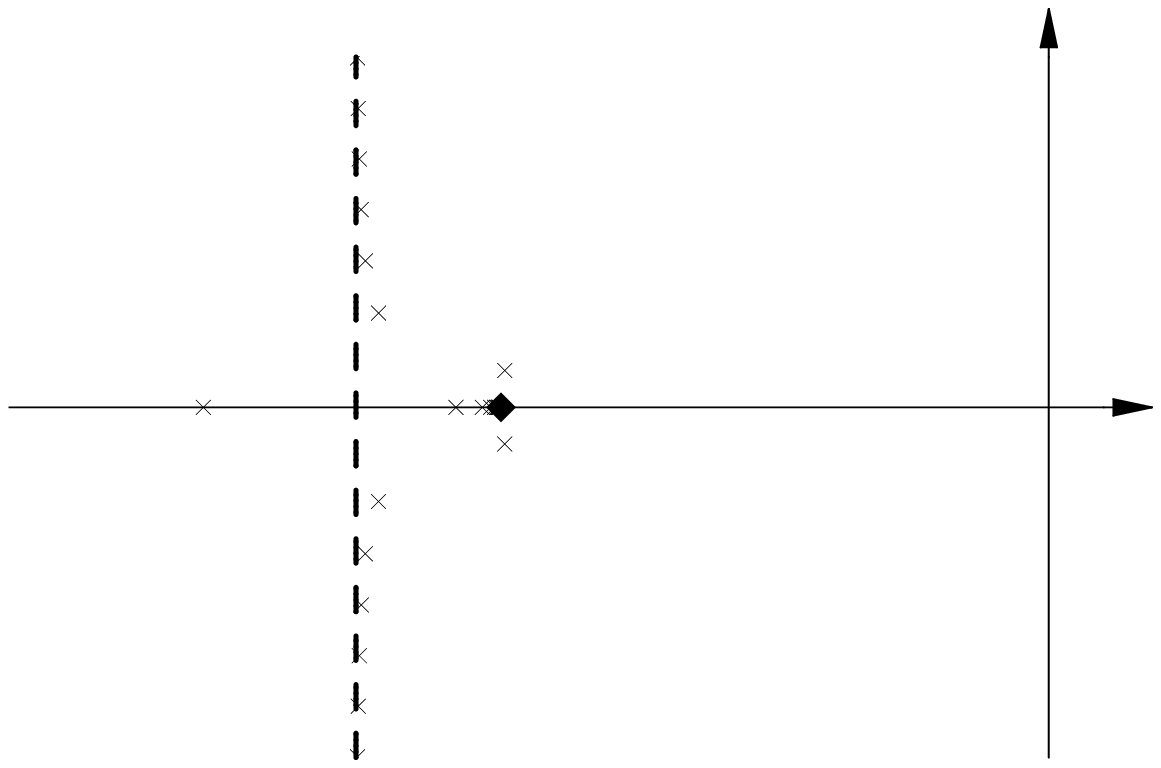}
                \caption{$\sigma_{max}=\Re (\lambda_2^1)=-1/\beta$}
                \label{fig:fig4}
        \end{subfigure}
\caption{Plots of the eigenvalues of the operator $A$ (crosses) in the complex plane (in solid lines, the real and complex axes), showing different possibilities for $\sigma_{max}(A)$. In all of them, the dashed line represents $\Re (\lambda) =-\frac{1}{2}\left( \frac{1}{\alpha}-\frac{1}{\beta}\right)$, which is the limit of the real parts of the nonreal eigenvalues, and the point marked as a diamond is $-\frac{1}{\beta}$, which is the limit of the real ones.}\label{fig:sigma_max}
\end{figure}

\begin{prop}[Dominant spectrum]\label{dominant}
 Let $\mu_1$ be the lowest eigenvalue of $L$. To find the dominant spectrum of $A$ one has to solve the cubic characteristic equation \eqref{eq:char_eq} with $\mu_n=\mu_1$. Then,
\begin{enumerate}
    \item
    \begin{enumerate}[(a)]
    \item If the three solutions of this equation are real (including the case of multiple solutions), then the dominant spectrum of $A$ will be $\left\{-1/\beta\right\}$ and $\sigma_{max}=\sigma_{max}(A)=-1/\beta$.
 \item If the solutions have the form $\lambda^1_1\in\R$ and $\lambda^1_2=\overline{\lambda^1_3}\in\C\setminus\R$, then the dominant spectrum of $A$ will be either $\{\lambda_2^1, \lambda_3^1\}$ or $\left\{-1/\beta\right\}$ (which is not an eigenvalue, but the only point in the essential spectrum), or both, depending on which has the highest real part. This real part is then $\sigma_{max}=\sigma_{max}(A)$.
     \end{enumerate}

\item We also claim that all the possibilities can occur as it is shown in the next three significative cases: \begin{enumerate}[(a)]
\item In the case $1/3\le\alpha/\beta<1$ one is in the situation 1(b) above and the dominant spectrum is $\{\lambda_2^{1}, \lambda_3^{1}\}$.
\item If $0<\alpha/\beta<1/3$ and $\mu_1$ is large enough one is in the situation 1(b) above but the dominant spectrum is $\left\{-1/\beta\right\}$.
\item If $0<\alpha/\beta<1/3$ and $\alpha$ is sufficiently small, with fixed $\beta$ and $\mu_1$, one is in the situation 1(a) above and so the dominant spectrum is $\left\{-1/\beta\right\}$.
    \end{enumerate}
    \end{enumerate}
In Figure \ref{fig:sigma_max} we can see different examples where the previous situations are attained.

\end{prop}

When the dominant spectrum is $\left\{-1/\beta\right\}$, then there will be no oscillations in the dominant part of the solutions and these cases could be qualified as {\em overdamped}.

Now we proceed to prove the previous propositions. Most of these results can be found in \cite{Trig2013} but we include the proofs here as they are slight but important differences in some of them.

\begin{pf}[of Proposition \ref{descriptioneig}]
First of all, let $\mu_n,\phi_n$ be a fixed eigenvalue and eigenfunction of $(L,D(L))$. We can look for solutions of \eqref{eq} of the form $u(x,t)=z(t)\phi_n(x)$. In this case, $z(t)$ would be a solution of
  $$    \alpha z'''+z''+\beta\mu_n z'+\mu_n z=0.$$
  This equation has \eqref{eq:char_eq} as its characteristic equation. Hence, for each $\mu_n$ eigenvalue of $(L,D(L))$ there exist three solutions of \eqref{eq:char_eq}, that form the three sequences $\lambda_1^n, \lambda_2^n, \lambda_3^n$ of eigenvalues of $(A,D(A))$, as stated. If any of these solutions is nonreal, we can write it as $\lambda_j^n=a_n+i\,b_n$, with $b_n\neq 0$. Imposing this in \eqref{eq:char_eq} and considering separately the real and the imaginary parts of the equation, it is easy to see that $a_n,b_n$ satisfy \eqref{eq:partreal} and \eqref{eq:partim}.

  To see part 1, we need to see whether this three solutions are real or not. For that, we will simply apply Cardano's method to \eqref{eq:char_eq} for a fixed $\mu_n$. The first thing to do in this method is to apply the change of variable $\xi=\lambda+\frac{1}{3\alpha}$  to \eqref{eq:char_eq} normalized such that the highest degree coefficient is equal to one. We obtain
  $$\xi^3+\xi p+q=0$$
  with
  \begin{equation}\label{eq:pq}
    p = \dfrac{\mu_n \beta}{\alpha}-\dfrac{1}{3\alpha^2},\ \
    q = \dfrac{2}{27\alpha^3}-\dfrac{\mu_n \beta}{3\alpha^2}+\dfrac{\mu_n }{\alpha}
  \end{equation}

  Now, it only remains to look at the sign of the Cardano's discriminant, that we can think as a function of $\mu_n$:
  \begin{equation}\label{eq:disc}
    d(\mu_n)= 4p^3+27 q^2 =\dfrac{\mu_n }{\alpha^2}\wwd
  \end{equation}
  with
  \begin{equation*}\label{eq:wdisc}
    \wwd = \dfrac{4 \beta^3}{\alpha}\mu_n^2 +\left(27-\left(\dfrac{\beta}{\alpha}\right)^2-18\left(\dfrac{\beta}{\alpha}\right) \right) \mu_n + \dfrac{4}{\alpha^2}.
  \end{equation*}
  According to Cardano's method:
  \begin{enumerate}[i)]
    \item if $d(\mu_n)>0$, the cubic polynomial has one real root and two complex conjugates
    \item if $d(\mu_n)<0$, the cubic polynomial has three different real roots
    \item if $d(\mu_n)=0$, all of the roots are real, with some of them being multiple.
  \end{enumerate}
  Observe that, as $\mu_n>0$, the sign of (\ref{eq:disc}) is determined by the sign of the second degree polynomial $\wwd$. The roots of this quadratic equation are the constants $m_1,m_2$ defined in \eqref{m1m2} with the constants $C_1,C_2$ defined in \eqref{eq:AB}. It can be seen that $C_2=0$ if and only of $\frac{\beta}{\alpha}=1$ or $\frac{\beta}{\alpha}=9$.
  So, depending on whether $m_1, m_2$ are real or not and positive (to coincide with a value of some $\mu_n$), we will have these different possibilities for the sign of the Cardano's discriminant (\ref{eq:disc}):

  \begin{enumerate}[(a)]

    \item Suppose $\frac{1}{9}<\frac{\alpha}{\beta}<1$. Then, one can see that $C_2<0$, which means that $\wwd$ will have a constant sign, which is positive. By Cardano's method, this concludes that \eqref{eq:char_eq} has one real root and the other two are complex conjugates for all $n\in\mathbb{N}$.

    \item Suppose $0<\frac{\alpha}{\beta}<\frac{1}{9}$. Then, one can see that $C_2>0$ and $C_1<0$. This means that $0<m_1<m_2$ are two different positive real roots for $\wwd$. This allows us to know the sign of $d(\mu_n)$:
        \begin{enumerate}[i)]
          \item $d(\mu_n)>0$ if $\mu_n\in(0,m_1)\cup(m_2,\infty)$. So, \eqref{eq:char_eq} has one real root and two complex conjugates if $n$ is such that $\mu_n\in(0,m_1)\cup(m_2,\infty)$ (we have an infinite number of them).
          \item $d(\mu_n)<0$ if $\mu_n\in(m_1,m_2)$. So, \eqref{eq:char_eq} has three different real roots for any possible $n$ such that $\mu_n\in(m_1,m_2)$ (if we have any of them, they will only be a finite number).
          \item $d(\mu_n)=0$ if there exist $n_1$ or $n_2$ such that $\mu_{n_1}=m_1$ or $\mu_{n_2}=m_2$. In this case, $\lambda_2^{n_1}=\lambda_3^{n_1}$ or $\lambda_2^{n_2}=\lambda_3^{n_2}$ is a real eigenvalue with algebraic multiplicity equal to two (and the other one is a simple real eigenvalue).
        \end{enumerate}
        Observe that the last two cases were not considered in \cite{Trig2013}.

    \item Finally, suppose that $\frac{\alpha}{\beta}=\frac{1}{9}$. In this case, in (\ref{eq:AB}) we have $C_2=0$ and $C_1=-216$, so $m_1=m_2=\frac{3}{\beta^2}>0$. This means that $d(\mu_n)>0$ for all $n\in\mathbb{N}$ except for, maybe, the case in which there exists $\mu_{n_1}=m_1=m_2$, for which the Cardano's discriminant is zero. But also observe that in this situation we also have $p=q=0$ (see (\ref{eq:pq})). Hence, by Cardano's method, \eqref{eq:char_eq} has one real root and the other two are complex conjugates for all $n\in\mathbb{N}$ except for, maybe, this certain $n_1$, for which $\lambda^{n_1}_1=\lambda^{n_1}_2=\lambda^{n_1}_3$ is a triple real root. If this triple root exists, then a simple computation allows to see that it is equal to $-\frac{1}{3\alpha}$, which is the same as $-\frac{3}{\beta}$. This case was not either considered in \cite{Trig2013}.
  \end{enumerate}

\vspace{0.1cm}

This proves part 1 of Proposition \ref{descriptioneig}. Let us now prove part 2. First, to prove \eqref{eq:fita_lambda}, we write as in \cite{Trig2013} $f(\lambda)= f_1(\lambda)+f_2(\lambda)$ with $f_1(\lambda)= \alpha \lambda^3+\lambda^2$ and $f_2(\lambda)= \beta\mu_n \lambda+\mu_n$, and we recall that we are interested in the real solutions of $f(\lambda)=0$. It is easy to see that for $\lambda\le -1/\alpha$ one has $f_1(\lambda)\le 0$ and $f_2(\lambda)<0$. Also, $f_1(\lambda)\ge 0$ for $\lambda\ge -1/\beta$ with $f_1(-1/\beta)>0$ and $f_2(\lambda)>0$ for $\lambda>-1/\beta$. So, all real roots of $f(\lambda)=0$ must be in $-1/\alpha<\lambda<-1/\beta$, as claimed. This is the same argument as that of \cite {Trig2013}, but we note that also holds in the case of three real eigenvalues.

Suppose now that $\lambda$ is a nonreal root of $f(\lambda)=0$. Then, for these values of $\alpha,\beta$ and $\mu_n$ the sign of the Cardano discriminant of \eqref{eq:char_eq} defined in \eqref{eq:disc} is positive. It is easy to see that the Cardano discriminant of \eqref{eq:partreal} is $d(\mu_n)/64$, with $d(\mu_n)$ defined in \eqref{eq:disc}. So, both discriminants have the same sign for the same values of the parameters, that is, positive in this case. Then, \eqref{eq:partreal} will also have a single real root, that will be precisely the real part of $\lambda$ we want to bound. We define $g(a)=8\alpha a^3+8 a^2+2 a\left( \frac{1}{\alpha}+\beta\mu_n\right)+\mu_n\left(\frac{\beta}{\alpha}-1\right)$, and we write \eqref{eq:partreal} as $g(Re(\lambda))=0$. We see that $g(-1/2(1/\alpha-1/\beta))<0$ and $g(0)>0$, and applying Bolzano's Theorem we conclude that \eqref{eq:fita_part re} holds.

To prove claim (b) of part 2 we consider now $\mu$ as a continuous variable in the open set $(0,m_1)\cup(m_2,\infty)$ or $(0,\infty)$ depending on $m_1$ and $m_2$ being real or not. According to what has been said above, in each of these open sets the equation $g(a)=0$ has a single and simple real solution $a=a_0(\mu)$. Deriving implicitly in $g(a_0(\mu))=0$ we obtain
$$\dfrac{da_0(\mu)}{d\mu}=\dfrac{-(2\beta a_0(\mu)+\frac{\beta}{\alpha}-1)}{g'(a_0(\mu))}.$$
Since we know that $a_0(\mu)>-1/2(1/\alpha-1/\beta)$, the numerator is negative. Since the coefficient of the cubic term of $g(a)$ is positive and $g(a)$ has only one real root, the derivative of $g$ at this root must be positive.

To conclude the proof of 2(b) we have still to prove that $a_0(\mu)$ decreases when $\mu$ jumps from $\mu=m_1$ to $\mu=m_2$. In more strict words, we want to show that $$\lim_{\mu\to m_1^-}a_0(\mu)\ge\lim_{\mu\to m_2^+}a_0(\mu).$$

Recall that these two values of $\mu$ are precisely the values for which the cubic equation $f(\lambda)=0$ has a double real root, and for $m_1<\mu<m_2$ the equation $f(\lambda)=0$ has three simple real roots, that we can order and call $\lambda_1(\mu)<\lambda_2(\mu)<\lambda_3(\mu)$. When $\mu\to m_1^+$, two of these roots collide, and become precisely $a_0(m_1)$, and the same happens as $\mu\to m_2^-$, when two of them collide and become $a_0(m_2)$. We do not know for the moment which of the three roots collide in each case, but it is clear that $\lambda_2(\mu)$ will be involved in the two collisions, since a collision between $\lambda_1(\mu)$ and $\lambda_3(\mu)$ is not possible without involving $\lambda_2(\mu)$. So, we conclude that $\lambda_2(\mu)\to a_0(m_1)$ as $\mu\to m_1^+$ and $\lambda_2(\mu)\to a_0(m_2)$ as $\mu\to m_2^-$.

So, our claim will be proved if we show that $\lambda_2'(\mu)<0$ for $m_1<\mu<m_2$. The central root of a cubic equation $f(\lambda)=0$ with three real simple roots and a positive coefficient of the cubic term is precisely the unique root that satisfies $f'(\lambda)<0$. Then, we can derive implicitly with respect to $\mu$ in the equation $f(\lambda_2(\mu))=0$ and obtain
$$\dfrac{d\lambda_2(\mu)}{d\mu}=\dfrac{-\beta\lambda_2(\mu)-1}{f'(\lambda_2(\mu))}.$$
The numerator is positive because of the upper bound in \eqref{eq:fita_lambda}, and the denominator is negative because of what we just said. This finishes the proof of part 2(b).

Finally, the proof of part 2(c) can be found in \cite{Trig2013} and part 2(d) is a straightforward computation.

\end{pf}

Let us now prove Proposition \ref{dominant}.

\begin{pf}[of Proposition \ref{dominant}]
To prove part 1(a) we observe that because of Proposition \ref{descriptioneig} all the real eigenvalues of $A$ satisfy $\lambda<-1/\beta$ and they accumulate at $-1/\beta$. To deal with the nonreal eigenvalues we observe that, under the hypotheses of 1(a), $m_1,m_2$ must be real and $m_1\le\mu_1\le m_2$. This, together with part 1(a) of Proposition \ref{descriptioneig}, implies that the values of $\mu_n$ that will give nonreal roots of \eqref{eq:char_eq} will satisfy $\mu_n>m_2$. Then, following the proof of Proposition \ref{descriptioneig} part 2(b), and with the same notation, $a_0(\mu_n)\le a_0(m_2)$. Even if $\mu=m_2$ is not an eigenvalue of $L$, the number $\lambda=a_0(m_2)$ will be a real (multiple) root of of \eqref{eq:char_eq} with $m_2$ in the place of $\mu_n$, so the bound \eqref{eq:fita_lambda} holds for this $\lambda$, $a_0(m_2)<-1/\beta$ and 1(a) is proved.

To prove 1(b) we have just to observe that, as we said, all the real eigenvalues satisfy $\lambda<-1/\beta$ and accumulate at this point, and, because of Proposition \ref{descriptioneig} part 2(b), the real parts of the nonreal ones are bounded above by $Re(\lambda_2^1)=Re(\lambda_3^1)$.

To prove 2(a) we observe, as we said in Proposition \ref{descriptioneig} part 2(d), that if $\alpha/\beta\ge 1/3$  then the vertical line $Re(z)=-1/2(1/\alpha- 1/\beta)$, where the nonreal eigenvalues accumulate (from its right hand side), lies at the right of the point $z=-1/\beta$ which is larger than all the real eigenvalues, so the real eigenvalues or their limit cannot be dominant.

Let us prove 2(b). Since $\alpha/\beta<1/3$, the point $z=-1/\beta$ lies at the right of the vertical line mentioned above. Since we know, by (\ref{eq:lim_Rel2}), that the function $a_0(\mu)$ as $\mu\to\infty$ tends to $-1/2(1/\alpha- 1/\beta)$, it is clear that if $\mu_1$ is large enough then $a_0(\mu_n)<-1/\beta-\varepsilon$ for all $n$ and some $\varepsilon>0$.

To prove part 2(c) we look at the expression of $d(\mu_1)$ as in \eqref{eq:disc} and observe that $d(\mu_1)<0$ if $\alpha$ is small enough, so we are in the situation 1(a) and the dominant spectrum will be $\{-1/\beta\}$.

\end{pf}

\section{A new scalar product, the normal property and decay of solutions}\label{sec:normal}

Let $\phi_n,\mu_n$ be the eigenfunctions and eigenvalues of $L$: $L \phi_n =\mu_n\phi_n$, in ascending order ($0<\mu_1\le\mu_2\le\dots\to\infty$) and with the collection $\{\phi_n\}$ being orthonormal in $H$. We consider the associated decompositions of the spaces in \eqref{espaisTrig} $\cH_i=\bigoplus_{n=1}^\infty E_n$, for $i=1,\dots 4$, where the $E_n$ are the three-dimensional spaces spanned by  $\{ (\phi_n,0,0), (0,\phi_n,0),(0,0,\phi_n)\}$. Observe that the spaces $E_n$ do not depend on $i$ and that for all $i$ they are orthogonal to each other with the natural scalar products of the spaces $\cH_i$.

These natural scalar products in the cases $i=1$ and $i=3$ when restricted to the spaces $E_n$ (and expressed in the previously given basis) are defined by the matrices

\begin{equation}\label{matrius}
O_{n,1} =\left(
            \begin{array}{ccc}
              \mu_n & 0 & 0 \\
              0 & \mu_n & 0 \\
              0 & 0 & 1 \\
            \end{array}
          \right),
\ \  O_{n,3} =\left(
            \begin{array}{ccc}
              \mu_n^2 & 0 & 0 \\
              0 & \mu_n & 0 \\
              0 & 0 & 1 \\
            \end{array}
          \right).
  \end{equation}
  We will focus only on $\cH_1$ and $\cH_3$ since the spaces $\cH_2$ and $\cH_4$ can be related to the previous ones by the natural isometry

  \begin{equation*}\label{isometry}
\left(
            \begin{array}{ccc}
              L^{1/2} & 0 & 0 \\
              0 & L^{1/2} & 0 \\
              0 & 0 & L^{1/2} \\
            \end{array}
          \right).
  \end{equation*}

Observe also that the scalar product \eqref{eq:prodMR} considered in \cite{MR2013}, is different, but equivalent, to $O_{n,1}$, as it can be written as

  \begin{equation*}\label{eq:MR}
    \left(
            \begin{array}{ccc}
              0 & 1 & 0 \\
              1 & \alpha & 0 \\
              0 & 1 & \alpha \\
            \end{array}
          \right)^T
          \left(
            \begin{array}{ccc}
              a^2\alpha(\beta-\alpha)\mu_n & 0 & 0 \\
              0 & a^2\mu_n & 0 \\
              0 & 0 & 1 \\
            \end{array}
          \right)
          \left(
            \begin{array}{ccc}
              0 & 1 & 0 \\
              1 & \alpha & 0 \\
              0 & 1 & \alpha \\
            \end{array}
          \right).
  \end{equation*}

  The idea is to define the new scalar product in each of the spaces $E_n$ by a symmetric real matrix $G_{n,i}$ expressed in the basis $\{ (\phi_n,0,0), (0,\phi_n,0),(0,0,\phi_n)\}$ in such a way that the eigenfunctions of $A$, that we call $\Psi^{n,i}_1,\Psi^{n,i}_2$ and $\Psi^{n,i}_3$, become orthonormal once normalized with the natural norm given by $O_{n,i}$. We understand that each $\Psi^{n,i}_j$ is the eigenfunction that corresponds to $\lambda^n_j$, once expressed in the previous basis. Hence, $\Psi^{n,i}_j=(1,\lambda^n_j,(\lambda^n_j)^2)^T/c^{n,i}_j$ where the $c^{n,i}_j$ is the normalizing constant in the usual norm depending on the space of \eqref{espaisTrig}. Then, we define the matrices

  \begin{equation}\label{Gni}
    C_{n,i}=\text{col}(\Psi^{n,i}_1,\Psi^{n,i}_2,\Psi^{n,i}_3) \textrm{ and }G_{n,i}=(\overline{C_{n,i}^{-1}})^TC_{n,i}^{-1}.
  \end{equation}
  When the matrix $C_{n,i}$ has the previous form it is easy to see that $G_{n,i}$ is a real, symmetric and positive definite matrix.

  The equivalence between the natural and the new norms is based on the following result.

  \begin{lem}\label{lem:quadrats}
  For $n$ sufficiently large, all $(x,y,z)\in\R^3$ (or $\C^3$) and all $i=1,\dots 4$ there exist numbers $M,m>0$ (independents of $n$) such that
  \begin{equation}\label{equivalents}
    m||(x,y,z)||_{O_{n,i}}\le ||(x,y,z)||_{G_{n,i}}\le M||(x,y,z)||_{O_{n,i}}
  \end{equation}
  \end{lem}

  \begin{pf}
  According to Proposition \ref{descriptioneig}, for sufficiently large $n$, the operator $A$ restricted to $E_n$ has eigenvalues $\lambda_1^n\in\R$ and $\lambda_2^n=a_n+ib_n=\overline{\lambda_3^n}\in\C\setminus\R$. By the same proposition, part 2 (c), the limits of $\lambda_1^n$, $a_n$ and $b_n$ are given by \eqref{eq:lim_l1}, \eqref{eq:lim_Rel2} and \eqref{eq:lim_Iml2}. Also, one can easily compute
  \begin{equation}\label{cs1}
  \begin{array}{l}
  c^{n,1}_1=\sqrt{ 1+\frac{1}{\beta^2} } \sqrt{\mu_n} +o(\sqrt{\mu_n} )\\
  c^{n,1}_2=c^{n,1}_3= \sqrt{\frac{\beta}{\alpha} + \frac{\beta^2}{\alpha^2} }\,\mu_n+o(\mu_n)\\
  \end{array}
\end{equation}
and
\begin{equation}\label{cs3}
  \begin{array}{l}
  c^{n,3}_1= \mu_n+o(\mu_n)\\
  c^{n,3}_2=c^{n,3}_3=\sqrt{1+\frac{\beta}{\alpha} +\frac{\beta^2}{\alpha^2}}\, \mu_n+o(\mu_n).\\
  \end{array}
\end{equation}
Then, one can compute the elements of the matrices $G_{n,1}$ and $G_{n,3}$. This is not a short calculation, and the use of an algebraic manipulator can be helpful. The results, up to leading orders as $\mu_n\to\infty$, are

\begin{equation}\label{Gn1}
  G_{n,1} =
   \left(
            \begin{array}{ccc}
               \frac{2\alpha\beta^2+3\alpha+\beta}{2\alpha\beta^2} \mu_n+ o(\mu_n) &  \frac{\alpha+\beta}{2\alpha\beta}\mu_n + o(\mu_n) &  \frac{4\alpha^2\beta^2+\beta^2+3\alpha^2}{4\alpha\beta^3} + o(1)  \\ \\
               \frac{\alpha+\beta}{2\alpha\beta}\mu_n + o(\mu_n) & \frac{\alpha+\beta}{2\alpha}\mu_n + o(\mu_n) & \frac{(\alpha+\beta)^2}{4\alpha\beta^2} + o(1) \\ \\
              \frac{4\alpha^2\beta^2+\beta^2+3\alpha^2}{4\alpha\beta^3} + o(1) & \frac{(\alpha+\beta)^2}{4\alpha\beta^2} + o(1)  & \frac{\alpha+\beta}{2\beta}+o(1) \\
            \end{array}
          \right)
\end{equation}

  and

\begin{equation}\label{Gn3}
  G_{n,3} =
   \left(
            \begin{array}{ccc}
               \mu_n^2 + o(\mu_n^2) & \frac{3\alpha\beta-\alpha^2+\beta^2}{2\alpha\beta^2} \mu_n +o(\mu_n) & \frac{\alpha}{\beta} \mu_n + o(\mu_n)\\ \\
               \frac{3\alpha\beta-\alpha^2+\beta^2}{2\alpha\beta^2} \mu_n +o(\mu_n) & \frac{\alpha^2+\alpha\beta+\beta^2}{2\alpha\beta} \mu_n + o(\mu_n) &
               \frac{6\alpha^2\beta-3\alpha^3+2\alpha\beta^2+\beta^3}{4\alpha\beta^3}+ o(1) \\ \\
              \frac{\alpha}{\beta} \mu_n + o(\mu_n) & \frac{6\alpha^2\beta-3\alpha^3+2\alpha\beta^2+\beta^3}{4\alpha\beta^3} + o(1) & \frac{3\alpha^2+\alpha\beta+\beta^2}{2\beta^2} + o(1) \\
            \end{array}
          \right).
\end{equation}

  Observe that all the leading terms of $G_{n,3}$ are positive because $\alpha<\beta$. Let us now first prove \eqref{equivalents} for $G_{n,1}$. Observe that, intuitively, this result will be true as both norms have the same diagonal terms (in asymptotic order) and the other ones are of lower order so they will be controlled by the diagonal ones. To prove that in a rigorous way, consider $(x,y,z)\in\mathbb{R}^3$ (or $\C^3$) and consider $n$ large enough such that all the terms of \eqref{Gn1} are positive (that is, large enough such that $o(\mu_n)$ and $o(1)$ do not affect the sign of the coefficients of the leading terms in $G_{n,1}$). Also, we are going to use two inequalities. First:

  \begin{equation}\label{quadrats_ineq}
  -\frac{1}{2}\left( c^2 a^2+\frac{1}{c^2}b^2\right) \leq ab \leq \frac{1}{2}\left( c^2 a^2+\frac{1}{c^2}b^2\right)
  \end{equation}
  which is true for any $a,b\in\R$ and $c>0$. And secondly, if $\{s_n\},\{r_n\}$ are positive real sequences such that $\lim_{n\to\infty} (s_n/r_n)=C>0$, it is easy to see that there exist $m,M>0$ such that
  \begin{equation}\label{lema5}
  m r_n \leq s_n \leq M r_n
  \end{equation}
  if $n$ is sufficiently large.

  Now, we start with the lower inequality of \eqref{equivalents}. From \eqref{Gn1} and using the left hand side inequality of \eqref{quadrats_ineq}

  \begin{equation}\label{lower}
  \begin{array}{c}
      (x,y,z) \,G_{n,1}\,(x,y,z)^T  \geq \left( \frac{2\alpha\beta^2+3\alpha+\beta}{2\alpha\beta^2} -c_1^2\,\frac{\alpha+\beta}{2\alpha\beta}+o(1) \right)\mu_n x^2   \\
        + \left( \frac{\alpha+\beta}{2\alpha} -\frac{1}{c_1^2}\frac{\alpha+\beta}{2\alpha\beta} + o(1) \right)\mu_n y^2 +
        \left( \frac{\alpha+\beta}{2\beta} - \frac{1}{c_2^2} \frac{4\alpha^2\beta^2+\beta^2+3\alpha^2}{4\alpha\beta^3} -  \frac{1}{c_3^2} \frac{(\alpha+\beta)^2}{4\alpha\beta^2} + o(1) \right) z^2
  \end{array}
  \end{equation}

  for certain $c_1,c_2,c_3>0$ such that the previous coefficients are positive for $n$ sufficiently large. Observe this is possible just taking $c_1>0$ and with $1/\beta < c_1^2 < \frac{2\alpha\beta^2+3\alpha+\beta}{\beta(\alpha+\beta)}$, and $c_2,c_3>0$ and large enough. Observe also that we can choose these constants independently of $n$.

  The same idea applies to prove the upper inequality of \eqref{equivalents}. From \eqref{Gn1} and using the right hand side inequality of \eqref{quadrats_ineq},

  \begin{equation}\label{upper}
  \begin{array}{c}
      (x,y,z) \,G_{n,1}\,(x,y,z)^T  \leq \left( \frac{2\alpha\beta^2+3\alpha+\beta}{2\alpha\beta^2} + c_4^2\,\frac{\alpha+\beta}{2\alpha\beta}+o(1) \right)\mu_n x^2   \\
        + \left( \frac{\alpha+\beta}{2\alpha} +\frac{1}{c_4^2}\frac{\alpha+\beta}{2\alpha\beta} + o(1) \right)\mu_n y^2 +
        \left( \frac{\alpha+\beta}{2\beta} + \frac{1}{c_5^2} \frac{4\alpha^2\beta^2+\beta^2+3\alpha^2}{4\alpha\beta^3} +  \frac{1}{c_6^2} \frac{(\alpha+\beta)^2}{4\alpha\beta^2} + o(1) \right) z^2
  \end{array}
  \end{equation}
  for certain $c_4,c_5,c_6>0$ such that the previous coefficients are positive for $n$ sufficiently large. In this case this is achieved simply taking $c_4=c_5=c_6=1$.

  Finally, with this choice of the constants $c_i$, \eqref{lema5} holds for \eqref{lower} and \eqref{upper}. Hence, there exist $m_1,M_1>0$ such that
  $$m_1 \left( \mu_n x^2 +\mu_n y^2 + z^2\right)\le ||(x,y,z)||^2_{G_{n,1}}\le M_1 \left( \mu_n x^2 +\mu_n y^2 + z^2\right)$$
  if $n$ is sufficiently large, which proves this lemma for $G_{n,1}$.

  The proof of this result for $G_{n,3}$ follows the same idea, but with some slight differences that we will point out. Again, consider $(x,y,z)\in\mathbb{R}^3$ (or $\C^3$) and $n$ large enough such that all the terms of \eqref{Gn3} are positive (that is, large enough such that $o(\mu_n^2)$, $o(\mu_n)$ and $o(1)$ do not affect the sign of the coefficients of the leading terms in $G_{n,3}$). We start with the lower inequality of \eqref{equivalents}. From \eqref{Gn3} and using the left hand side inequality of \eqref{quadrats_ineq},

  \begin{equation}\label{lowerGn3}
  \begin{array}{c}
      (x,y,z) \,G_{n,3}\,(x,y,z)^T  \geq \left( \mu_n^2-c_1^2 \frac{3\alpha\beta-\alpha^2+\beta^2}{2\alpha\beta^2}\mu_n -c_2^2 \frac{\alpha}{\beta} \mu_n + o(\mu_n^2) \right) x^2   \\
        + \left( \frac{\alpha^2+\alpha\beta+\beta^2}{2\alpha\beta}\mu_n -\frac{1}{c_1^2} \frac{3\alpha\beta-\alpha^2+\beta^2}{2\alpha\beta^2}\mu_n
         - c_3^2\frac{6\alpha^2\beta-3\alpha^3+2\alpha\beta^2+\beta^3}{4\alpha\beta^3}+ o(\mu_n) \right) y^2 \\ +
        \left( \frac{3\alpha^2+\alpha\beta+\beta^2}{2\beta^2} -\frac{1}{c_2^2}(\frac{\alpha}{\beta}\mu_n + o(\mu_n))-\frac{1}{c_3^2}\frac{6\alpha^2\beta-3\alpha^3+2\alpha\beta^2+\beta^3}{4\alpha\beta^3}   + o(1) \right) z^2
  \end{array}
  \end{equation}
  for other $c_1,c_2,c_3>0$ such that the previous coefficients are positive for $n$ sufficiently large and of order $O(\mu_n^2)$, $O(\mu_n)$ and $O(1)$, respectively. For this to be true, we will need to choose $c_2=c_2(\mu_n)$. It suffices to choose $c_1>0$, independent of $n$ and such that $c_1^2>\frac{3\alpha\beta-\alpha^2+\beta^2}{\beta(\alpha^2+\alpha\beta+\beta^2)}$, $c_2=\widetilde{c_2}\sqrt{\mu_n}$ with $\widetilde{c_2}>0$, independent of $n$ and such that $\frac{2\alpha\beta}{3\alpha^2+\alpha\beta+\beta^2} <\widetilde{c_2}^2<\frac{\beta}{\alpha} $, and $c_3>0$, independent of $n$ and such that $c_3^2> \frac{6\alpha^2\beta-3\alpha^3+2\alpha\beta^2+\beta^3}{4\alpha\beta^2}\left(\frac{3\alpha^2+\alpha\beta+\beta^2}{2\beta} -\frac{\alpha}{\widetilde{c_2}^2} \right)^{-1}$.

  For the upper inequality, from \eqref{Gn3} and using the right hand side inequality of \eqref{quadrats_ineq} one gets
  \begin{equation}\label{upperGn3}
  \begin{array}{c}
      (x,y,z) \,G_{n,3}\,(x,y,z)^T  \leq \left( \mu_n^2+c_4^2 \frac{3\alpha\beta-\alpha^2+\beta^2}{2\alpha\beta^2}\mu_n +c_5^2 \frac{\alpha}{\beta}\mu_n + o(\mu_n^2) \right) x^2   \\
        + \left( \frac{\alpha^2+\alpha\beta+\beta^2}{2\alpha\beta}\mu_n +\frac{1}{c_4^2} \frac{3\alpha\beta-\alpha^2+\beta^2}{2\alpha\beta^2}\mu_n
         + c_6^2\frac{6\alpha^2\beta-3\alpha^3+2\alpha\beta^2+\beta^3}{4\alpha\beta^3}+ o(\mu_n) \right) y^2 \\ +
        \left( \frac{3\alpha^2+\alpha\beta+\beta^2}{2\beta^2} +\frac{1}{c_5^2}(\frac{\alpha}{\beta}\mu_n + o(\mu_n))+\frac{1}{c_6^2}\frac{6\alpha^2\beta-3\alpha^3+2\alpha\beta^2+\beta^3}{4\alpha\beta^3}   + o(1) \right) z^2
  \end{array}
  \end{equation}
  for other $c_4,c_5(\mu_n),c_6>0$ such that the previous coefficients are positive for $n$ sufficiently large and of the right order. In this case this is achieved simply taking $c_4=c_6=1$ and $c_5=\sqrt{\mu_n}$.

  So, as in the case of $G_{n,1}$, with the previous choice of the new constants $c_i$, \eqref{lema5} also holds for \eqref{lowerGn3} and \eqref{upperGn3}. Hence, there exist $m_3,M_3>0$ such that
  $$m_3 \left( \mu_n^2 x^2 +\mu_n y^2 + z^2\right)\le ||(x,y,z)||^2_{G_{n,3}}\le M_3 \left( \mu_n^2 x^2 +\mu_n y^2 + z^2\right)$$
  if $n$ is sufficiently large, which proves the present lemma also for $G_{n,3}$.
  \end{pf}

Let us now prove Theorem \ref{thm1}.

\begin{pf} [of Theorem \ref{thm1}]
When $\mu_n\ne m_1,m_2$, then in each of the $A$-invariant three-dimensional subspaces $E_n$ defined in the beginning of this Section there are three different eigenvalues of $A$ and one can consider the scalar product given by the real symmetric matrices $G_{n,i}$ defined in \eqref{Gni}. We can then extend the definition of the scalar product to the whole of $\cH_i=\bigoplus_{n=1}^\infty E_n$, by a block-diagonal procedure, $G_i=\text{diag}(G_1,G_2,\dots)$.

With the scalar product so defined, the subspaces $E_n$ are orthogonal to each other, so the whole set of eigenfunctions
  \begin{equation}\label{eigenfunctions}
    F^i=\{\Psi^{n,i}_1,\Psi^{n,i}_2,\Psi^{n,i}_3; \ n=1,2\dots\}
  \end{equation}
becomes orthonormal. The operator $A$ diagonalizes in this basis, its adjoint $A^*$ is given by just its conjugate matrix, and so $A$ and $A^*$ commute and $A$ is a normal operator.

To see that this new scalar product gives a norm that is equivalent to the old natural norm, we use the Lemma \ref{lem:quadrats} above in $E_{n_0+1}\oplus E_{n_0+2}\oplus\dots$ for $n_0$ large enough and use in $E_{1}\oplus E_{2}\oplus\dots E_{n_0}$ that in finite dimensions all norms are equivalent.

To finish the reasoning we have still to prove that the family $F^i$ of eigenfunctions is complete in each $\cH_i$. Suppose that $U'=(u',v',w') \in  \cH_i$ is a nonzero vector that is orthogonal to the whole family $F^i$, and we will arrive to a contradiction. If $U'$ is nonzero, then at least one of its three components will be a nonzero element of $H$. Because of that, it will have at least a nonzero component in the basis $\{\phi_n\}$, suppose for $n=n'$, so it will have a nonzero projection in $E_{n'}$, namely $U_{n'}$. Since the projection $U'\mapsto U_{n'}$ is orthogonal in all of the cases $i=1,2,3$ and $4$, we see that $\langle U',\Psi^{n',i}_j\rangle_{G_i}=\langle U_{n'},\Psi^{n',i}_j\rangle_{G_i}$, and this cannot be zero for all $j=1,2,3$ if $U_{n'}$ is nonzero.

To prove the last part of the Theorem let us suppose $\mu_{n_1}=m_1<m_2$ (the cases $m_1<m_2=\mu_{n_2}$ or $\mu_{n_1}=m_1=m_2$ are similar). In this case the characteristic equation \eqref{eq:char_eq} has one double real root $\lambda_2^{n_1}=\lambda_3^{n_1}$ and a different simple real root $\lambda_1^{n_1}$. The restriction of $A$, as it appears in \eqref{eq:evolution}, to the invariant subspace $E_{n_1}$ expressed in the basis $\{ (\phi_{n_1},0,0), (0,\phi_{n_1},0),(0,0,\phi_{n_1})\}$ will have the form $$
    \left(
      \begin{array}{ccc}
        0 & 1 & 0 \\
        0 & 0 & 1 \\
        -\frac{\mu_{n_1}}{\alpha} & -\frac{\beta\mu_{n_1}}{\alpha} & -\frac{1}{\alpha} \\
      \end{array}
    \right).$$

It is easy to see that this matrix has $\lambda_2^{n_1}$ as an eigenvalue of geometric multiplicity one but algebraic multiplicity two. This is a property that will hold independently of the scalar product considered. And it is well known that this is impossible for normal operators, that have the property that geometric and algebraic multiplicities of eigenvalues always coincide.
  \end{pf}

Let us now proceed with the proof of Theorem \ref{thm2}.

\begin{pf}[of Theorem \ref{thm2}]
The proof of this theorem is the same in all the spaces given in \eqref{espaisTrig}. Hence, our notation will no distinguish among them and we will not include the superindex $i$, which distinguishes among the spaces.
\begin{enumerate}[i)]
\item For the parameter values that make $A$ a normal operator in the suitable new scalar product $G$ given in Theorem \ref{thm1}, it has been shown that there exists an orthonormal and complete set of eigenfunctions $\{ \Psi_j^n\}$, with $A \Psi_j^n = \lambda_j^n \Psi_j^n$, $j=1,2,3$, $n=1,\ldots,\infty$. If $U(0)=\sum_{n,j} d_j^n \Psi_j^n$, then $U(t)=\sum_{n,j} d_j^n e^{\lambda_j^n t} \Psi_j^n$ and, because of the orthonormality of the eigenfunctions,
    $$ \|U(t)\|_G^2=\sum_{n,j} |d_j^n|^2 e^{2 Re(\lambda_j^n) t} \leq \sum_{n,j} |d_j^n|^2 e^{2 \sigma_{max} t} \textrm{   } (t>0).$$

\item We have seen in Proposition \ref{dominant} that $\sigma_{max}$ is either $\Re (\lambda_2^{1})$ or $-1/\beta$. In the first case, the solution $U(t)=e^{\lambda^{1}_2 t}\, \Psi_2^{1}$ itself has the optimal decay rate. In the second case, if $\sigma_{max}=-1/\beta$, the sequence $\lambda_1^n$ tends to $-1/\beta$ from the left (see Proposition \ref{descriptioneig}, parts 2(a) and 2(c)) and the corresponding solutions $U_n(t)=e^{\lambda_1^{n}t} \Psi_1^{n}$ have decay rates $\lambda_1^n$, which can be taken as close as we want to $-1/\beta$.

\item The idea of the proof of this part is that when there are non-semisimple eigenvalues they cannot be dominant.
    To proceed in this way, among the sequence $0<\mu_1\leq \mu_2\leq \cdots \mu_n\leq\cdots\to\infty$ of eigenvalues of $L$  we distinguish the finite set $S$ of those that coincide either with $m_1$ or $m_2$ defined in \eqref{m1m2} (see Proposition \ref{descriptioneig}), parts 1(b) and 1(c)) and accordingly decompose $H=H_0 \oplus H_1$ and
    $L=\left(
         \begin{array}{cc}
           L_0 & 0 \\
           0 & L_1 \\
         \end{array}
       \right)$
    in such a way that $\sigma(L_0) =\sigma(L)\setminus S$ and $\sigma(L_1)=S$ ($H_1$ is finite dimensional). We make the same corresponding decomposition in each of the spaces given in \eqref{espaisTrig}, $\cH_i=\cH_i^0\oplus \cH_i^1$, and the operator
    $A=\left(
         \begin{array}{cc}
           A_0 & 0 \\
           0 & A_1 \\
         \end{array}
       \right)$.
    Observe that $A_0$ is in the situation described in the first part of Theorem $\ref{thm1}$ and $A_1$ is a finite dimensional operator, with all its eigenvalues being real, and some of them being multiple.
    Hence, according to Theorem \ref{thm1}, we can define a new scalar product $\langle \cdot,\cdot\rangle_{G_0}$ in $\cH_i^0$ in which we can obtain the optimal decay inequality of part i) above
    $$ \|e^{A_0 t}\|_{G_0} \leq e^{\sigma_{max}(A_0)\,t}\ \ \textrm{ for }t\geq 0.$$
    On the other hand, as $\cH_i^1$ is finite-dimensional and according to a well-know result of Linear Algebra, for each $\varepsilon>0$ we can define a new scalar product $\langle \cdot,\cdot\rangle_{G_{1,\varepsilon}}$ in $\cH_i^1$ such that
    $$ \|e^{A_1 t}\|_{G_{1,\varepsilon}} \leq e^{\left(\sigma_{max}(A_1)+\varepsilon\right)\,t}\ \ \textrm{ for }t\geq 0.$$
    As it is deduced from Proposition \ref{descriptioneig} part 2(a), $\sigma_{max}(A_1) <-1/\beta$. So, we can choose $\varepsilon>0$ such that $\sigma_{max}(A_1)+\varepsilon<-1/\beta\leq \sigma_{max}(A_0)$.
    Finally, we define the scalar product $G_i'$ in $\cH_i$ as the orthogonal extension of $G_0$ and $G_1$. It is equivalent to the natural scalar product of each $\cH_i$ because it is so when restricted to each of $\cH_i^0$ and $\cH_i^1$. And the optimal decay rate result follows in the $G_i'$ norm because the dominant part of the spectrum is in $\sigma(A_0)$ and the optimality is true for $G_0$ because of part ii) above.

\end{enumerate}

\end{pf}

\end{document}